\documentclass[natbib,11pt]{article}

\RequirePackage[OT1]{fontenc}

\usepackage{amssymb}
\usepackage{theorem}
\usepackage{xspace}
\usepackage[final]{graphicx}
\usepackage{subfigure}
\usepackage{xspace}
\usepackage{hhline}

 \usepackage{amsmath}
\usepackage{bm}
\usepackage{natbib}
\usepackage{bbm}
\usepackage[a4paper]{geometry}

\usepackage{geometry}

\theoremstyle{plain}

\usepackage{natbib}   
\usepackage{array}
 \usepackage{subfigure}
\usepackage{epsfig}

\def\efo{\textrm{\mathversion{bold}$\mathbf{\phi^0}$\mathversion{normal}}}

\def\eS{\textrm{\mathversion{bold}$\mathbf{\Sigma}$\mathversion{normal}}} 

\def\eE{I\!\!E}
\def\eP{I\!\!P}
\def\e1{1\!\!1}

\def\eg{\overset{.}{\textbf{g}}}
\def\egg{\overset{..}{\textbf{g}}}

\def\ev{\mathbf{{v}}}
\def\xx{\textrm{\mathversion{bold}$\mathbf{x}$\mathversion{normal}}}
\def\XX{\textrm{\mathversion{bold}$\mathbf{X}$\mathversion{normal}}}
\def\uu{\textrm{\mathversion{bold}$\mathbf{u}$\mathversion{normal}}}

\def\pth#1{\left(#1\right)}
\def\acc#1{\left\{#1\right\}}
\def\cro#1{\left[#1\right]}

\newtheorem{theorem}{Theorem}[section]

\newtheorem{lemma}{Lemma}[section]
\newtheorem{proposition}{Proposition}[section]
\newtheorem{remark}{Remark}[section]

\def\ef{\textrm{\mathversion{bold}$\mathbf{\phi}$\mathversion{normal}}}  
\def\ee1{\textrm{\mathversion{bold}$\mathbf{\varepsilon}$\mathversion{normal}}}  
\def\eth{\textrm{\mathversion{bold}$\mathbf{\theta}$\mathversion{normal}}}

\newcommand{\N}{\mathbb{N}}
\newcommand{\R}{\mathbb{R}}
\newcommand{\Z}{\mathbb{Z}}

\newcommand{\Var}{\mathbb{V}\mbox{ar}\,}

\def\argmin{\mathop{\mathrm{arg\,min}}} 
\def\hh{ \hspace*{0.5cm}}
 
\begin{document}

\title {{\bf Estimation in a change-point nonlinear quantile model}}


\author{GABRIELA CIUPERCA  \footnote{Universit\'e Lyon 1,  Institut Camille Jordan, 
Bat.  Braconnier, 43, blvd du 11 novembre 1918, 
F - 69622 Villeurbanne Cedex, France,  E-mail: {\it Gabriela.Ciuperca@univ-lyon1.fr},}\\
Universit\'e Lyon 1,  UMR 5208, Institut Camille Jordan,  France
}
\maketitle



{\textbf{ABSTRACT. }}
{ {\normalsize
\noindent 
\textit{This paper considers a nonlinear quantile model with change-points. The quantile estimation method, which as a particular case includes median model, is more robust with respect to other traditional methods when model errors contain outliers. Under relatively weak assumptions, the convergence rate and asymptotic distribution of change-point and of regression parameter estimators are obtained.  Numerical study by Monte Carlo simulations shows the performance of the proposed method for nonlinear model with change-points. }
}}\\



 



\section{Introduction}
\label{sec1:1}
\hh Classically, for linear or nonlinear models, the errors are supposed with mean zero and   bounded variance. In this case, model parameters are estimated generally by least squares (LS) method. If these conditions are not satisfied or if model contains outliers, then the LS estimators of the model parameters can have a large error. A very interesting and robust alternative method was proposed by \cite{Koenker:Bassett:78} by the introduction of the quantile method. A particular case of this method is that of least absolute deviation (LAD). For a complete overview on quantile method, we refer the reader to book of \cite{Koenker:05}. Properties of a nonlinear quantile model are studied also in the papers \cite{Chen:Tran:Lin:04}, \cite{Choi:Kim:Lee:05}, \cite{Oberhofer:Haupt:13}. \\
\hh On the other hand, in applications, it is possible that we have not one but several models, the localization where model changes being unknown. We obtain which is called as  a change-point model. The purpose of this paper is to study the properties of this type of model  estimated by quantile method,  when between two consecutive change-points the model is nonlinear. For this study we need to known the asymptotic behaviour of the objective function. \\
\hh To our knowledge, most previous studies of change-point models have focused on linear models. On this subject, we can mention the following papers: \cite{Bai:98} for LAD method,  \cite{Bai:Perron:98} for LS method, \cite{Koul:Qian:02} for maximum likelihood method, \cite{KQS} for M-estimation method. For quantile method, \cite{Oka:Qu:11} estimate the change-points location and the coefficient parameters of each phase, \cite{Furno:12} realize a Lagrange multiplier test for detecting the structural breaks. For change-point nonlinear model, because of difficulties caused by the nonlinearity, literature is less rich: \cite{Boldea:Hall:13} use LS method to estimate and test the number of breaks. In  \cite{Ciuperca:09}, the M-estimation method is used to estimate a multiphase nonlinear model with random design and changes in the model due to some (unknown) values in design.  A general criterion is proposed in \cite{Ciuperca:11a} to determine the change-point number. If changes in the model occur in time, the LAD estimation method was studied by  \cite{Ciuperca:11}. \\
\hh Present paper generalizes \cite{Ciuperca:11}, considering a method, for estimating and for choosing the change-point number criterion, based on the quantile framework. This is because, often in practice, especially in the case of change-point models, the quantile index $\tau$ of errors is not 1/2.\\
\hh We note the important fact that, in a multiple change-point model, the change-point estimation could affects the estimator properties. Moreover, it is difficult to study, theoretically but also numerically, a change-point model since it depends of two parameter types: the regression  parameters and the change-points. \\ 
\hh The plan of this paper is as follows. In Section \ref{sec2:1} we first introduce some notations and assumptions. Next, we study the asymptotic behaviour of the objective function. In Section \ref{sec3:1}, we define and study the quantile estimator in a nonlinear model with change-points. Convergence rate and asymptotic distributions of the estimators are obtained. Finally, in Section \ref{sec4:1}, simulation results illustrate the performance of the quantile method for change-point nonlinear model. In Appendix Bernstein's inequality is recalled.    
\section{Quantile regression without change}
\label{sec2:1}
\hh In this section we  study the asymptotic behaviour of the quantile process. \\

Let us consider the  following regression model
\begin{equation}
\label{eq1}
Y_i=g(\XX_i,\ef)+\varepsilon_i, \qquad i=1, \cdots , n,
\end{equation}
where the regression function $g: \Upsilon \times \Gamma \rightarrow \R$, with $\ef \in \Gamma \subseteq \R^p$,  $\xx \in \Upsilon$, $\Upsilon \subseteq \R^d$ is known up to the parameter $\ef$. We suppose that the set $\Gamma $ is compact.\\
\hh For a fixed quantile index $\tau \in (0,1)$, the $\tau$th conditional quantile regression of $Y$, given $\xx$, is $g(\xx,\ef)+F^{-1}(\tau)$, with $F^{-1}(\tau)$  the $\tau$th quantile ($F^{-1}$ is the inverse of the distribution function  $F$) of error $\varepsilon$. We suppose that $F(0)=\tau$. \\
In the case when the model contains intercept, noted $\phi_0$, the regression function has the form   $g(\xx,\ef)=\phi_0+g_0(\xx,\ef_1)$. Then, the following parameter vector is considered $\ef(\tau)= (\phi_0(\tau), \ef_1)^t$, with $\phi_0(\tau) = \phi_0+F^{-1}(\tau)$. Thus, when model contains an intercept $\phi_0$, we estimate first $\phi_0(\tau)$ from which we then have the estimation of $\phi_0$. For linear models, in literature,  the presence or not of intercept intervenes in the result proofs (see for example \cite{Oka:Qu:11}).\\
\hh The estimator of $\ef(\tau)$ by  Koenker and Bassett technique (see \cite{Koenker:Bassett:78}) is called  \textit{quantile regression}. We suppose $\tau$ fixed, therefore, for simplicity reasons,  we will note $\ef(\tau)$ by $\ef$. Contrary to the papers where linear models are studied  (see for example \cite{Oka:Qu:11}) when  $\phi_0\neq 0$, in this paper we shall consider   simultaneously  the two cases presented above.

\subsection{Assumptions and notations}
 \hh In this subsection we give assumptions and notations needed in this paper. \\
\hh For simplicity, we suppose that the regressors $\XX_i$ are non random, although the results will, typically hold for random $\XX_i$'s independent of the $\varepsilon_i$'s and if $\XX_i$ independent of $\XX_j$ for $i \neq j$. \\
For the  model (\ref{eq1}), we consider  $\ef^0$  the true value  (unknown) of $\ef$, with $\efo$ an inner point of the compact $\Gamma$. \\
For a fixed quantile index  $\tau \in (0,1)$, consider the check function $\rho:\R \rightarrow \R$ given by $$\rho_\tau(u)  = u [\tau - \e1_{u \leq 0}]$$ and  random variable $$D_i(\tau) \equiv \tau- \e1_{\varepsilon_i \leq 0}.$$ Since $F(0)=\tau$ we have $\eE[D_i(\tau)]=0$. \\
\hh The quantile  estimator of  parameter $\ef$ is defined by 
\begin{equation}
\label{ett}
\hat \ef_n^{(\tau)} \equiv \argmin_{\ef \in \Gamma} \sum^n_{i=1} \rho_\tau ( Y_i - g(\XX_i,\ef)).
\end{equation}
 Its consistency, with $n^{-1/2}$  convergence rate   and asymptotic normality of estimator $\hat \ef_n^{(\tau)}$  have been proved in previous papers (see for  example \cite{Koenker:05}).   \\
 
 For reading convenience, all throughout the paper,  let us consider the following notation, for all $\ef \in \Gamma$ and $i=1, \cdots, n$, 
  $$h_i(\ef)\equiv g(\XX_i,\ef)- g(\XX_i,\ef^0).$$
 In order to study the quantile model, let us consider the following two random processes:
\begin{equation}
\label{eq2}
\left\{
\begin{array}{l}
\displaystyle{ W_n(\tau;\ef,\efo) \equiv -\sum^n_{i=1} D_i(\tau) h_i(\ef)} ,\\
\displaystyle{ Z_n(\ef,\efo) \equiv \sum^n_{i=1} \int^{h_i(\ef)}_0
 (\e1_{\varepsilon_i \leq s}-\e1_{\varepsilon_i \leq 0})ds .}
 \end{array} 
\right.
\end{equation}
Obviously $\eE[W_n(\tau;\ef,\efo)]=0$. \\
\hh For each sample $i \in \{1, \cdots , n \}$, in order to study always the quantile model,  consider the following difference 
\[
G_i^{(\tau)}(\ef, \ef^0) \equiv \rho_\tau (\varepsilon_i-g(\XX_i,\ef)+ g(\XX_i,\ef^0) )  - \rho_\tau(\varepsilon_i),  
\]
 from which  us define the random   process 
$${\cal G}^{(\tau)}_n(\ef, \ef^0) \equiv \sum^n_{i=1}G_i^{(\tau)}(\ef, \ef^0),$$
which is in fact the objective function for finding the quantile estimator $\hat \ef^{(\tau)}_n$ of model (\ref{eq1}). \\
  Using identity of  \cite{Knight:98}, for  any real nonzero number $r$, we have that 
$
\rho_\tau(r-s)- \rho_\tau(r)=s[ \e1_{r <0} - \tau] +\int^{s}_0 [\e1_{r \leq t} - \e1_{r \leq 0}]dt$. 
Then, the process  ${\cal G}^{(\tau)}_n(\ef, \ef^0)$ can also write as
\begin{equation}
\label{eq3}
{\cal G}^{(\tau)}_n(\ef, \ef^0) = W_n(\tau;\ef,\efo) +Z_n(\ef,\efo).
\end{equation}
In the case of a nonlinear model,  function ${\cal G}^{(\tau)}_n(\ef, \ef^0) $ is not  convex in $\ef$  which means that the study of estimator  $\hat \ef_n^{(\tau)}$ and of function  ${\cal G}^{(\tau)}_n(\ef, \ef^0) $ will be different than for a linear model, based on the convexity.\\

All throughout the paper, C denotes a positive generic constant which may take different values in different formula or even in different parts of the same formula. 
 All vectors are column and  $\textbf{v}^t$ denotes the transposed of $\textbf{v}$. All vectors and matrices are in bold.  Concerning the used norms, for a m-vector $\ev=(v_1, \cdots,v_m)$, let us denote by   $\|\ev\|_1= \sum   _{j=1} ^ m |v_j|$  and $\|\ev\|_2=(\sum   _{j=1} ^ m v^2_j)^{1/2}$. For a matrix $ \textbf{M}=(a_{ij})_{\substack{1\leqslant i \leqslant m_1\\1 \leqslant j \leqslant m_2}}$, we denote by $ \| \textbf{M}\|_1 = \max   _ {j=1,\cdots, m_2} (\sum  _{i=1}^{m_1} |a_{ij}|)$, the subordinate norm to the vector norm $\| .\|_1$.\\
 
We now state the assumptions on the errors and on the regression function.\\
\hh The errors $(\varepsilon_i)_{1 \leq i \leq n}$ are supposed independent identically distributed (i.i.d.) random variables.  We denote by  $f$ the density and by  $F$ the  distribution function of  $\varepsilon$. \\
\textbf{(A1)} There exists two constants $c_0>0$ and  $\delta_0>0$ such that for all $|x| \leq \delta_0 $, we have 
\[
\min(F(|x|)-F(0),F(0)-F(-|x|)) \geq c_0 |x|.
\]
Contrary to the classic assumptions for a nonlinear regression model, we do not impose the condition that the mean of errors $\varepsilon_i$ is zero or that their variance is bounded. \\
\hh The regression function $g(\xx,\ef)$ is supposed twice differentiable in $\ef$ and continuous on $\Upsilon$. In the following, for $\xx \in \Upsilon$ and $\ef \in \Gamma$ we use notation $\eg(\xx,\ef) \equiv \partial g(\xx, \ef)/ \partial \ef$ and $\egg(\xx,\ef) \equiv \partial^2 g(\xx, \ef)/ \partial \ef^2$. Moreover, for the function $g$, following assumptions are considered:\\
\textbf{(A2)} For all $\xx \in \Upsilon$, function $\eg(\xx,\ef)$ is bounded in every  $\eta$-neighbourhood  of $\efo$, when  $\eta \rightarrow 0$.\\
\textbf{(A3)} There exists $c_1>0$ such that  $n^{-1} \sum^n_{i=1} \sup_{\ef \in \Gamma} \| \eg(\XX_i,\ef) \|_2 \leq c_1 < \infty$.\\
\textbf{(A4)} There exist two positive constants $c_2,c_3 >0$ and natural $n_0$ such that for all $\ef_1, \ef_2 \in \Gamma$ and $n \geq n_0$:
$
c_2 \| \ef_1 - \ef_2 \|_2 \leq \left( n^{-1} \sum^n_{i=1} [g(\XX_i,\ef_1) -g(\XX_i,\ef_2)]^2 \right)^{1/2} \leq c_3 \| \ef_1 - \ef_2 \|_2.
$
Moreover, we have
 $ n^{-1}   \sum^n_{i=1} \eg(\XX_i,\efo) \eg^t(\XX_i,\efo)$ converges, as $n \rightarrow \infty$,  to a positive definite matrix. Furthermore, $\max_{1 \leq i \leq n} n^{-1/2}  \| \eg(\XX_i,\efo) \|_2 \rightarrow 0$.\\
\textbf{(A5)} For all $\xx \in \Upsilon$, $\ef \in {\cal V}_\eta(\efo)$, with $\eta\rightarrow 0$, we have that   $\|\egg(\xx,\ef)\|_1$ is bounded. \\
For certain results, stronger assumptions are necessary:\\
\textbf{(A6)} For all $\ef \in \Gamma$, $\xx \in \Upsilon$, we have that  $\|\eg(\xx,\ef)\|_2$ is bounded.\\
\textbf{(A7)} For all $\ef \in \Gamma$, $\xx \in \Upsilon$, we have that  $\|\egg(\xx,\ef)\|_1$ is bounded. \\

We wish to emphasise the fact that, with respect to the particular case $\tau=1/2$, considered  in   \cite{Ciuperca:11}, we consider here  multidimensional regressors  $\XX_i$.  \\
\hh Assumption (A2) means that for every $\eta \rightarrow 0$, the function $\|\eg(\xx,\ef)\|_2$ is bounded for all $\xx \in \Upsilon $ and for all $\ef \in {\cal V}_\eta$, with 
$${\cal V}_\eta \equiv \{ \ef \in \Gamma ;  \| \ef- \efo \|_2 \leq \eta \}.$$
\hh Assumptions (A1), (A4) are needed that the objective function has an unique minimum at $\efo$ and for  convergence and asymptotic normality  of the quantile estimator (see Koenker (2005), page 124). Obviously that, assumption  (A7) implies (A5). We have also that  (A6) implies (A2), (A3) and third condition of (A4).\\
Note that  assumption (A3) is the same as in paper \cite{Ciuperca:11}, for a median nonlinear model and (A1) is supposition (C4) of  Oberhofer and Haupt (2014)'s paper, for a nonlinear quantile regression with weakly dependent errors. As noted in the last paper, assumption (A1) is stronger that the usual assumption in literature: $f(x)$ exists in a neighbourhood of $x=0$ and $f(0) \geq c_0 >0$.
 
\subsection{Asymptotic behaviour of objective function ${\cal G}_n^{(\tau)}(\ef,\efo)$}
\hh In order to study the main part of this paper devoted to a change-point model estimated by quantile framework, in this subsection  the asymptotic behaviour of process ${\cal G}_n^{(\tau)}(\ef,\efo)$ is studied. \\
\hh Recall that under assumptions (A1) and (A4), the quantile estimator $\hat \ef^{(\tau)}_n$, given by (\ref{ett}), is weakly $n^{-1/2}$-consistent (see Koenker (2005) or Oberhofer (1982)). \\

\begin{remark}
By elementary calculations we show that, for all $\ef \in \Gamma$:
\begin{itemize}
\item if $h_i(\ef) \geq 0$ we have 
\[\int_0^{h_i(\ef)}\e1_{\varepsilon_i <t} dt=h_i(\ef) \e1_{\varepsilon_i <0}+[h_i(\ef) -\varepsilon_i] \e1_{0 \leq \varepsilon_i \leq h_i(\ef)}+0 \cdot \e1_{\varepsilon_i >h_i(\ef)},\]
\item if $h_i(\ef) < 0$ we have
\[
\int_0^{h_i(\ef)}\e1_{\varepsilon_i <t} dt=0 \cdot \e1_{\varepsilon_i >0} - \varepsilon_i \e1_{h_i(\ef) \leq \varepsilon_i \leq 0} +h_i(\ef) \e1_{\varepsilon_i <0}.\]
\end{itemize}
These imply that
\begin{itemize}
\item if $h_i(\ef) \geq 0$ we have\\
\[ \int_0^{h_i(\ef)}[\e1_{\varepsilon_i <t} -\e1_{\varepsilon_i<0}]dt 
 =[h_i(\ef)-\varepsilon_i] \e1_{0 \leq \varepsilon_i \leq h_i(\ef)},\]
\item if $h_i(\ef) < 0$ we have 
\[\int_0^{h_i(\ef)}[\e1_{\varepsilon_i <t} -\e1_{\varepsilon_i<0}]dt =[-\varepsilon_i -h_i(\ef) ] \e1_{h_i(\ef) \leq \varepsilon_i \leq 0}.\]
\end{itemize}
Thus, for $Z_n(\ef,\efo)$ defined by (\ref{eq2}), we have with probability one that, $Z_n(\ef,\efo) \geq 0$ for all $\ef \in \Gamma$.
\end{remark}

A consequence of this remark, taking into account that $\eE[W_n(\tau;\ef,\efo)] =0$ for all $\ef \in \Gamma$, is that 
\begin{equation}
\label{egn}
\eE[{\cal G}_n^{(\tau)}(\ef,\efo)] \geq 0, \qquad  \textrm{ for all } \ef \in \Gamma.
\end{equation}

  We will now study the asymptotic behaviour of the objective function ${\cal G}_n^{(\tau)}(\ef,\efo)$. In this purpose, for a  bounded deterministic sequence $c_n$, let us consider the following parameter set: $$\Omega_{c_n} \equiv \{ \ef \in  \Gamma;   \| \ef-\efo\|_2 \leq c_n \}.$$
Emphasize that for the following  Proposition,  claim  (i), the sequence  $(c_n)$ is a  constant $c$. Thus, we denote the set  $\Omega_{c_n}$ by $\Omega_c$. The proof idea is the same as in \cite{Bai:98}, Lemma 4, only now the nonlinearity of $g(\xx,\ef)$ and the quantile regression intervene significantly.
\begin{proposition}
\label{Lemma2.4_adapt} 
Let us consider a deterministic positive sequence  $(a_n)$ such that $a_n \rightarrow \infty $,   as $n \rightarrow \infty$. \\
(i) Under assumption (A6), if sequence  $(a_n)$ satisfies  in addition the conditions  $n^{-1} a_n =O(1)$, $n^{-1}a_n^2/\log n \rightarrow \infty$ as $n \rightarrow \infty$, and the parameters belong to the set $\Omega_c =\{\ef \in \Gamma ; \| \ef-\efo \|_2 \leq c \}$,   we have that, for all $\epsilon >0$, 
there exists a constant $C>0$ and a natural number $n_\epsilon \in \N$ such that for all $ n \geq n_\epsilon$,
\[
 \eP \cro{\sup_{\ef \in \Omega_{c}} \left|\frac{1}{a_n} \cro{ {\cal G}_n^{(\tau)}(\ef,\efo)- \eE[{\cal G}_n^{(\tau)}(\ef,\efo)]}  \right| >\epsilon } 
 \leq \exp (- C \epsilon^2 n^{-1} a^2_n  ).
\]
(ii) Under assumptions (A2), (A3), if we have furthermore an another sequence $(c_n)$ such that $c_n \rightarrow 0$, $n^{-1}c^{-1}_n a_n =O(1)$, $a_n^2/(n c_n^2 \log n) \rightarrow \infty$ as $n \rightarrow \infty$ and the parameters belong to the set $\Omega_{c_n}$, 
 we have that, for all $\epsilon >0$, 
there exists a constant $C>0$ and a natural number $n_\epsilon \in \N$ such that for all $ n \geq n_\epsilon$,
\[
 \eP \cro{\sup_{\ef \in \Omega_{c_n}} \left|\frac{1}{a_n} \cro{ {\cal G}_n^{(\tau)}(\ef,\efo)- \eE[{\cal G}_n^{(\tau)}(\ef,\efo)]}  \right| >\epsilon } 
 \leq \exp (- C \epsilon^2 n^{-1} c_n^{-2} a^2_n  ).
\]
\end{proposition}
\noindent {\bf Proof.}
 \textit{(i)} We decompose $\Omega_c$ in subsets, such that the diameter of each subset is less than $c n^{-1/2}$. Thus, we can write $\Omega_c =\cup^{C_p n^{p/2}}_{j=1} \omega_j^n$, where $\omega^n_j \equiv \{  \ef \in \Omega_c ; \| \ef - \ef_j \|_2 \leq c n^{-1/2}\}$, with  $ \ef_j \in \Omega_c $.
This decomposition depends on the dimension of the parameter  set $\Gamma$, then positive constant $C_p$($< \infty$) depends on $p$. Diameter $c n^{-1/2}$ of each cell $\omega^n_j$ was taken as that, for two points belonging to the same cell, the difference which occurs in the left hand side of  following relation (\ref{eth1}) converges to 0 as $n \rightarrow \infty$. Then, in order to study the behaviour of ${\cal G}_n^{(\tau)}(\ef,\efo)- \eE[{\cal G}_n^{(\tau)}(\ef,\efo)]$ for all $\ef \in \Gamma$,  we take one representative point, noted $\ef_j$, of $\omega^n_j$,  for each $j =1, \cdots, C_p n^{p/2}$. Then, in order to study the probability of the left hand side of claim \textit{(i)} we just study the probability of the left hand side of (\ref{eth2}).\\
Since $n^{-1} a_n^2 / \log n \rightarrow \infty$, then we have that $n^{-1/2} a_n \rightarrow \infty$ as $n \rightarrow \infty$. Thus, for each $\ef_1, \ef_2 \in \omega^n_j$, using assumption   $n^{-1} \sum^n_{i=1} \sup_{\ef \in \Gamma}$ $ \| \eg(\XX_i,\ef) \|_2 \leq c_1 < \infty$ of (A6), we have
  \begin{equation}
  \label{eth1}
  \begin{array}{l}
  a_n^{-1} |{\cal G}_n^{(\tau)}(\ef_1,\efo)- \eE[{\cal G}_n^{(\tau)}(\ef_1,\efo)]   - {\cal G}_n^{(\tau)}(\ef_2,\efo)+ \eE[{\cal G}_n^{(\tau)}(\ef_2,\efo)] |  \\
\qquad \qquad  \qquad  \leq C  n a_n^{-1} \| \ef_1 - \ef_2 \|_2  \leq C n^{1/2} a_n^{-1} \rightarrow 0.
   \end{array}
  \end{equation}

  For  $\ef_j \in \omega^n_j$, for any $j =1, \cdots, C_p n^{p/2}$, we have the following obvious inequality
\begin{equation}
\label{eth2}
\displaystyle{\eP\left[\sup_j a_n^{-1} |{\cal G}_n^{(\tau)}(\ef_j,\efo)- \eE[{\cal G}_n^{(\tau)}(\ef_j,\efo)] |  > \epsilon \right]  }  
\displaystyle{\leq  \sum^{C_p n^{p/2}}_{j=1} \eP \left[|{\cal G}_n^{(\tau)}(\ef_j,\efo)- \eE[{\cal G}_n^{(\tau)}(\ef_j,\efo)] |  > \epsilon a_n \right]. }
\end{equation}
\hh On the other hand, taking into account assumption (A6), we obtain that there exists a constant $C>0$ such that 
 $|G_i^{(\tau)}(\ef_j,\efo) - \eE[G_i^{(\tau)}(\ef_j,\efo)] |  \leq C \| \ef_j-\efo\|_2 \leq C c$, with probability 1.  \\
Moreover $\Var [{\cal G}_n^{(\tau)}(\ef_j,\efo)] \leq C n c^2$, with $C >0$.\\
Then, since $n^{-1} a_n =O(1)$,  $a^2_n/(n  \log n) \rightarrow \infty$   as $n \rightarrow \infty$, we can apply Bernstein's inequality  (\ref{ineq_Bernstein}) (see Appendix), for $\beta = C c$, $V=Cn c^2$, $s=1/2$, $z=a_n \epsilon$ and we obtain that 
\[
\begin{array}{ll}
\displaystyle{ \sum^{C_p n^{p/2}}_{j=1}\eP [|{\cal G}_n^{(\tau)}(\ef_j,\efo)-\eE[{\cal G}_n^{(\tau)}(\ef_j,\efo)] |  > \epsilon a_n]} & \leq 2 C_p n^{p/2} \exp(-C \epsilon^2 n^{-1} a^2_n) \\
& = 2 C_p   \exp(-C \epsilon^2 n^{-1} a^2_n+p/2 \log n) \rightarrow 0,
\end{array}
\]
 with $C>0$.\\
  The claim follows combining the last relation together with (\ref{eth1}) and (\ref{eth2}).\\
\hh  \textit{(ii)} In this  case, we decompose the set $\Omega_{c_n} =\cup^{C_p n^{p/2}}_{j=1} \omega_j^n$, with the subsets  \\ $\omega^n_j \equiv \{  \ef \in \Omega_{c_n} ; \| \ef - \ef_j \|_2 \leq  c_n  n^{-1/2}\}$, with $ \ef_j \in \Omega_{c_n} $.\\
Since $a^2_n / (n c^2_n \log n ) \rightarrow \infty$ we have that $n^{1/2} c_n/a_n \rightarrow 0$ as $n \rightarrow \infty$. Then, for $\ef_1, \ef_2 \in \omega^n_j$ we have by similar arguments as in (\ref{eth1}): $a_n^{-1} |{\cal G}_n^{(\tau)}(\ef_1,\efo)- \eE[{\cal G} _n^{(\tau)}(\ef_1,\efo)] - {\cal G}_n^{(\tau)}(\ef_2,\efo)+ \eE[{\cal G}_n^{(\tau)}(\ef_2,\efo)] | \leq C n^{1/2} c_n a_n^{-1} \rightarrow 0$, as $n \rightarrow \infty $.\\
  On the other hand, for $\ef_j \in \omega^n_j$, for any  $j =1, \cdots, C_p n^{p/2}$, using assumption (A2), (A3) we have: $|G_i^{(\tau)}(\ef_j,\efo) - \eE[G_i^{(\tau)}(\ef_j,\efo)] | \leq  C c_n$, with probability  1. The end of proof is similar to that of (i), using $\Var [{\cal G}_n^{(\tau)}(\ef_j,\efo)] \leq C n c^2_n$ and  applying   Bernstein's inequality  (\ref{ineq_Bernstein})  for $\beta = C c_n$, $V=Cn c^2_n$, $s=1/2$, $z=a_n \epsilon$.
 \hspace*{\fill}$\blacksquare$ \\
 
 \begin{remark} An example of sequence $(a_n)$ which satisfies the conditions of Proposition \ref{Lemma2.4_adapt}(i) is  $a_n=n /(\log n)^s$, with $s \geq 0$.
 \end{remark}

Applying Proposition \ref{Lemma2.4_adapt} and   Borel-Cantelli lemma we obtain that for all $\epsilon >0$, we have
\begin{equation}
\label{eq38_adapt}
 \limsup_{n \rightarrow \infty} \left(\sup_{ \ef \in \Omega_{c_n}} \frac{1}{a_n} \left| {\cal G}_n^{(\tau)}(\ef, \ef^0) - \eE[{\cal G}_n^{(\tau)}(\ef, \ef^0)]\right|\right)  \leq \epsilon,  \qquad a.s.
\end{equation}
which is the equivalent of Lemma A.2 of \cite{Oka:Qu:11}.\\

The following Proposition proves a general result that the infimum of objective function ${\cal G}_n^{(\tau)}(\ef, \ef^0)$ is strictely positive for $\ef$ outside a ball centred at $\efo$ with radius $w_n$ , when $w_n \rightarrow 0$.

\begin{proposition}
\label{etape4}
Suppose that assumptions  (A1), (A4)-(A6) are satisfied,  that density $f$ of $\varepsilon$ is differentiable in a neighbourhood of 0 and that  $f'(x) $ is bounded in this neighbourhood.\\
Let $(w_n)$ be a monotone deterministic sequence   converging to  0, such that $n w^2_n / \log n \rightarrow \infty$, as $n \rightarrow \infty$.\\
 Then, there exists $\epsilon_1 >0$ such that,
\[
 \liminf_{n \rightarrow \infty} \pth{ \inf_{\|\ef -  \ef^0\|_2  \geq  w_n} \frac{1}{n w^2_n}  {\cal G}_n^{(\tau)}(\ef, \ef^0) } > \epsilon_1, 
\]
with  probability 1.
\end{proposition}
\noindent {\bf Proof.}
Let us first consider a monotone positive sequence $(v_n)$ such that $v_n \rightarrow 0$,  as $n \rightarrow \infty$ and $v_n \geq w_n $ for $n$ large enough.\\
We now consider parameter   $\ef$ such that $\|\ef-\efo\|_2 =v_n$.  
For a  $p$-vector  $\uu$ in an open set of $\R^p$, with $\|\uu\|_2=1$, we have using (\ref{eq3}),
\[
\eE[{\cal G}^{(\tau)}_n (\efo+v_n \uu ,\efo)]=\eE[Z_n(\efo+ v_n \uu, \efo)]=\sum^n_{i=1}  \int^{ h_i(\efo+ v_n \uu)}_0 [F(s)-F(0)] ds.
\]
We make  the Taylor's expansion for the distribution function $F$  and we obtain
 \begin{equation}
\label{R1R2}
 \eE[{\cal G}^{(\tau)}_n (\efo+v_n \uu ,\efo)] 
   =  \sum^n_{i=1}  \int^{  h_i(\efo+ v_n \uu)}_0  s f(0) ds 
    + \sum^n_{i=1}  \int^{ h_i(\efo+ v_n \uu)}_0   \frac{s^2}{2} f'( \zeta_i s) ds  
 \equiv  R_1 +R_2,
\end{equation}
with  $0< \zeta_i <1$. \\
The Taylor's expansion up to order 1 of $g$ at $\ef =\efo$, the assumptions (A5), (A6) imply that  $R_1$ is equal to
\begin{equation}
\label{R1}
 v^2_n \frac{f(0)}{2}  \sum^n_{i=1} [ \uu^t \eg(\XX_i,\efo)+v_n  \uu^t \egg(\XX_i, \tilde \ef_i) \uu   ]^2  = v^2_n \frac{f(0)}{2}  \sum^n_{i=1} [ \uu^t \eg(\XX_i,\efo)]^2(1+o(1)),
\end{equation}
with $\tilde \ef_i=\efo+b_i v_n \uu$, $b_i \in [0,1]$, $i=1, \cdots, n$.\\
For $R_2$, by similar arguments and since $f'$ is bounded in the neighbourhood of  0, we have
\begin{equation}
\label{R2}
|R_2|  \leq  C \sum^n_{i=1} \int^{ |h_i(\efo+ v_n \uu)|}_0 \frac{s^2}{2} ds 
  =  C v^3_n \sum^n_{i=1} [ \uu^t \eg(\XX_i,\efo)]^3(1+o(1))
  =  o(R_1).
\end{equation}
Then, taking into account (\ref{R1R2}), (\ref{R1}) and (\ref{R2}) we obtain:

$$
\eE[{\cal G}^{(\tau)}_n (\efo+v_n \uu ,\efo)]= w^2_n \frac{f(0)}{2} \sum^n_{i=1} [ \uu^t \eg(\XX_i,\efo) ]^2(1+o(1)),$$ 
with $o(1)$ uniformly in $\uu$. 
Since $\| \uu \|_2=1$, we have then that
\[
\eE[{\cal G}^{(\tau)}_n (\efo+ v_n \uu,\efo) ] \geq n w^2_n \frac{f(0)}{2} \lambda_{\min, n}(1+o(1)),
\]
where $\lambda_{\min, n}$ is  the smallest eigenvalue of matrix  $n^{-1} \sum^n_{i=1} \eg(\XX_i,\efo) \eg^t(\XX_i,\efo)$. \\
On the other hand, assumption  (A4) implies  that there exists a $\lambda >0$ such that $\lambda_{\min, n}\rightarrow \lambda$, as $n \rightarrow \infty$. Then, since $n v^2_n \rightarrow \infty$, we have that for all  large enough  $n$:
\begin{equation}
\label{eq_etoile2}
\frac{1}{n v^2_n} \eE[{\cal G}^{(\tau)}_n (\efo+v_n \uu ,\efo) ] \geq  \lambda \frac{f(0)}{2}(1+o(1)).
\end{equation}
Under assumptions (A6), we take  $\epsilon =\lambda f(0)/4 $, $\Omega_{n} \equiv \{ \ef \in \Gamma; \| \ef-\efo\|_2=v_n \}$ and $a_n=n v^2_n$ for Proposition \ref{Lemma2.4_adapt}(ii).  Then,  relation (\ref{eq38_adapt}) follows. Moreover,  since  $(n v^2_n)^{-1} \Var [{\cal G}^{(\tau)}_n (\ef, \efo)] < \infty$, $n v^2_n / \log n \rightarrow \infty$, together with relation (\ref{eq_etoile2}), we have that
\begin{equation}
\label{vn}
\liminf_{n \rightarrow \infty} ( \inf_{\|\ef -  \ef^0\|_2 = v_n} (n v^2_n)^{-1} {\cal G}_n^{(\tau)}(\ef, \ef^0) ) > \eta,
\end{equation}
 with $\eta = \lambda  f(0)/4$. \\
 Let us now consider a monotone positive sequence $(v_n)$ such that $v_n$ (or no sub-sequence) does not converge to 0. Then,  there exists $\delta >0$ such that $v_n \geq \delta$ for any $n$ large enough.\\
 By simple algebraic computations we obtain that:
\[
\eE[G_i^{(\tau)}(\ef,\efo)]=\left\{
\begin{array}{l}
 \int^{-h_i(\ef)}_{0}[|h_i(\ef)| -x] dF(x),  \textrm{ if }  h_i(\ef) <0, \\
 \int_{-h_i(\ef)}^{0}[|h_i(\ef)| +x] dF(x), \textrm{ if }   h_i(\ef) \geq 0,
\end{array}
\right.
\]
 Then 
 \[
 \begin{array}{ll}
 \eE[G_i^{(\tau)}(\ef,\efo)] & \displaystyle{  \geq  \e1_{h_i({\boldsymbol{\phi}}) \geq 0}\int_{-\frac{h_i({\boldsymbol{\phi}})}{2}}^0 \left(|h_i({\boldsymbol{\phi}})|+x\right) dF(x) +\e1_{h_i({\boldsymbol{\phi}}) <0} \int^{-\frac{h_i({\boldsymbol{\phi}})}{2}}_0 (|h_i({\boldsymbol{\phi}})|-x) dF(x)}   \\ 
  & \displaystyle{ \geq   \e1_{h_i({\boldsymbol{\phi}}) \geq 0} \frac{|h_i({\boldsymbol{\phi}})|}{2} \int_{-\frac{h_i({\boldsymbol{\phi}})}{2}}^0 dF(x) +\e1_{h_i({\boldsymbol{\phi}}) <0} \frac{|h_i({\boldsymbol{\phi}})|}{2} \int^{-\frac{h_i({\boldsymbol{\phi}})}{2}}_0 dF(x)}   \\
   & \displaystyle{   =  \frac{|h_i({\boldsymbol{\phi}})|}{2} \left(\e1_{h_i({\boldsymbol{\phi}}) \geq 0} [F(0)-F(-\frac{h_i({\boldsymbol{\phi}})}{2})] +\e1_{h_i({\boldsymbol{\phi}}) <0}  [F(-\frac{h_i({\boldsymbol{\phi}})}{2})-F(0)]\right).}  
 \end{array}
\]
 
Taking into account the assumptions (A1) and (A4), for $\| \ef-\efo\|_2=v_n \geq \delta$, we have that there exits a constant $c>0$ such that:
 \begin{equation}
 \label{eq_etoile2_bis}
 \eE[n^{-1} {\cal G}^{(\tau)}_n(\ef,\efo) ] >c.
 \end{equation}
 Under assumption  (A6), taking  $\epsilon =c/2 $, $\Omega_{n} \equiv \{ \ef \in \Gamma; \| \ef-\efo\|_2=v_n \}$ and $a_n=n v^2_n$, for Proposition \ref{Lemma2.4_adapt}, we obtain the relation (\ref{eq38_adapt}). Then,  together with relation (\ref{eq_etoile2_bis}) we have that
\begin{equation}
\label{vn_bis}
\liminf_{n \rightarrow \infty} ( \inf_{\|\ef -  \ef^0\|_2 = v_n} (n v^2_n)^{-1} {\cal G}_n^{(\tau)}(\ef, \ef^0) ) > \eta,
\end{equation}
 with $\eta = \epsilon =c/2$. \\
 Since relations (\ref{vn}) and (\ref{vn_bis}) are valid for any positive sequence $(v_n)$, such that $v_n \geq w_n$, and since $(w_n)$ is monotonic, the Proposition follows if we consider  $\epsilon_1 = \eta$.

 \hspace*{\fill}$\blacksquare$ \\

 The following two lemma will be needed in the next section, where model contains change-points. The  change-points are the observations where model changes. In the next section, we will  estimate simultaneously these change-points but also the model parameters between two change-points. The following Lemma will be used to find the convergence rate of the  change-points estimators. 
 \begin{lemma}
\label{Lemma2.3_adapt}
For $1 \leq l < k \leq n$ such that $k-l \rightarrow \infty $, as $n \rightarrow \infty$, under assumptions (A2), (A7),  we have, for all $\alpha >1/2$,
$$
\sup_{1 \leq l < k \leq n} | \inf_{\ef} \sum^k_{i=l} G_i^{(\tau)}(\ef,\efo) | =O_{\eP}(n^\alpha).$$ 
\end{lemma}
\noindent {\bf Proof.}
Since for all  $r_1, r_2 \in \R$ we have that $|\rho_\tau(r_1) - \rho_\tau(r_2)| < |r_1-r_2|$, then 
$$| \rho(\varepsilon_i -h_i(\ef_1))-\rho(\varepsilon_i -h_i(\ef_2)) | \leq | h_i(\ef_1) - h_i(\ef_2)| .$$
 For $\|\ef_1- \ef_2\|_2 \leq C n^{-1/2}$, using the Taylor expansion up to order 1 for each of the two functions $h_i(\ef_1)$, $h_i(\ef_2)$ in respect to $\ef$, around  $\efo$, and using assumptions (A2), (A7),   we obtain that 
$
\sum^n_{i=1} \{ \rho(\varepsilon_i -h_i(\ef_1))-\rho(\varepsilon_i -h_i(\ef_2)) -\eE[\rho(\varepsilon_i -h_i(\ef_1))-\rho(\varepsilon_i -h_i(\ef_2))] \} \leq O_{\eP}(n^{1/2})$.\\
Since $\eE[G_i^{(\tau)}(\ef, \ef^0) ] \geq 0$ for all $\ef$, we have that
 $G_i^{(\tau)}(\efo, \ef^0) =0 \geq \inf_{\ef} \sum^k_{i=l}G_i^{(\tau)}(\ef, \ef^0) \geq \inf_{\ef} \sum^k_{i=l}\left( G_i^{(\tau)}(\ef, \ef^0)- \eE[ G_i^{(\tau)}(\ef, \ef^0) ]\right) $ .\\
  The rest of proof is similar  to that of   Lemma 2.3  of \cite{Ciuperca:14}.
  \hspace*{\fill}$\blacksquare$ \\

By the following Lemma we prove that the objective function ${\cal G}^{(\tau)}_n(\ef,\efo)$ given by (\ref{eq3}) varies little  when a  small portion of observations is ignored. 
\begin{lemma}
\label{Lemma L+}
Under assumptions (A2), (A5), for all parameter $\ef$ such that  $\| \ef-\efo \|_2 \leq n^{-1/2}$ and for $M \in \N$ arbitrary, we have, as $n \rightarrow \infty$, 
$$
\sup_{m\in [n-M,n]} \sup_{\ef, \|\ef-\efo\|_2 \leq n^{-1/2} } |{\cal G}^{(\tau)}_m(\ef,\efo)-{\cal G}^{(\tau)}_n(\ef,\efo) |=o_{\eP}(1).$$ 
\end{lemma}
\noindent {\bf Proof.}
Let us consider a natural number  $m$ such that $m\in [n-M,n]$. 
Since $|\rho_\tau(\varepsilon-r)-\rho_\tau(\varepsilon)| \leq |r|$, we have that
$|{\cal G}^{(\tau)}_m(\ef,\efo)-{\cal G}^{(\tau)}_n(\ef,\efo) |  \leq  \sum^n_{i=m+1} | g(\XX_i,\ef)-g(\XX_i,\efo)| \leq  \sum^n_{i=m+1} \|\ef - \efo \|_1 \| \eg(\XX_i,\efo) +\egg(\XX_i,\tilde \ef)(\ef- \efo)^t\|_1$, 
 with $\tilde \ef = \efo + a (\ef- \efo)$, $a \in [0,1]$. Since $\egg$ is bounded  in a neighbourhood of $\efo$ by assumption (A5),  applying also the Markov inequality and assumption (A2), the last sum is smaller than $n^{-1}O_{\eP}(1)=O_{\eP}(n^{-1})$.
\hspace*{\fill}$\blacksquare$ \\

\begin{lemma}
 \label{Lemma 5_LAD}
 Under the same assumptions as in Proposition \ref{etape4}, we have, for all $\delta \in (0,1)$:
 \[
 \sup_{[n \delta] \leq m \leq n} \left| \inf_{\ef, \| \ef-\efo\|_2 \leq n^{-1/2} }  {\cal G}^{(\tau)}_m(\ef,\efo) \right| =O_{\eP}(1).
 \]
  \end{lemma} 
 {\bf Proof.}\\
 Given the convergence rate  of the  quantile estimator, we have
 \[
 \sup_{[n \delta] \leq m \leq n}\left| \argmin_{\ef \in  \Gamma} {\cal G}^{(\tau)}_m(\ef,\efo) \right|=O_{\eP}(n^{-1/2}).
 \]
Then,  it is sufficient to prove that
 \begin{equation}
 \label{ibis}
 \sup_{1 \leq m \leq n} \sup_{\| \uu\|_2 \leq M}|{\cal G}^{(\tau)}_m(\efo+n^{-1/2} \uu,\efo)  |=O_{\eP}(1).
 \end{equation}
 On the other hand, by relation (\ref{eq2}), for any $\uu \in \R^p$ such that $\| \uu\|_2 \leq M$, we have  ${\cal G}^{(\tau)}_m(\efo+n^{-1/2} \uu,\efo) =W_m(\tau; \efo+n^{-1/2} \uu,\efo)+Z_m(\efo+n^{-1/2} \uu,\efo)$, with,   by definition: $Z_m(\efo+n^{-1/2} \uu,\efo)=\sum^m_{i=1} \int^{h_i(\efo+n^{-1/2} \uu)}_0 [\e1_{\varepsilon_i \leq s}-\e1_{\varepsilon_i \leq 0}] ds$ and
 $W_m(\tau; \efo+n^{-1/2} \uu,\efo)=-\sum^m_{i=1} D_i(\tau) h_i(\efo+n^{-1/2} \uu)$.\\
 We first study  $W_m$. Recall that $\eE[W_m(\tau; \efo+n^{-1/2} \uu,\efo)]=0$. On the other hand, taking the Taylor expansion up to order 1 of  $g$, we have:
 \[ W_m(\tau; \efo+n^{-1/2} \uu,\efo)= - \sum^m_{i=1}D_i \uu^t  [n^{-1/2} \eg(\XX_i,\efo)+n^{-1}  \egg(\XX_i,\tilde \ef_i) \uu]\]
 with $\tilde \ef_i=\efo+b_i n^{-1/2} \uu$, $b_i \in [0,1]$, $i=1, \cdots, n$.\\
By the central limit theorem, using assumptions (A4) and (A6),  we have:
 \begin{equation}
 \label{nDi}
 n^{-1/2} \sum^m_{i=1} D_i \uu^t \eg(\XX_i,\efo)=O_{\eP}(1), 
 \end{equation}
 uniformly in $m$ and $\uu$. Since $\| \egg(\XX_i,\tilde \ef_i)\|_1$ is bounded by assumption (A5), we obtain:
 \begin{equation}
 \label{nug}
 n^{-1} \sum^m_{i=1} \uu^t \egg(\XX_i,\tilde \ef_i) \uu =O_{\eP}(1),
 \end{equation}
  uniformly in $m$ and $\uu$. Thus, by (\ref{nDi}) and (\ref{nug}), we have 
  \begin{equation}
  \label{Wmm}
  W_m(\tau; \efo+n^{-1/2} \uu,\efo)=O_{\eP}(1),  
  \end{equation}
  uniformly in  $m$  and $\uu$.\\
   We study now $Z_m(\efo+n^{-1/2} \uu,\efo)$:
   \[
   \begin{array}{l}
   \eE[Z_m(\efo+n^{-1/2} \uu,\efo) ] 
   \displaystyle{ =   \sum^m_{i=1} \int^{ h_i(\efo+ n^{-1/2} \uu) }_0 [F(s)-F(0)] ds   } \qquad \\
   \displaystyle{ \qquad   =  \sum^m_{i=1}  \int^{ h_i(\efo+ n^{-1/2} \uu)}_0  s f(0) ds+ \sum^n_{i=1}  \int^{ h_i(\efo+ n^{-1/2} \uu)}_0   \frac{s^2}{2} f'( \zeta_i s) ds  }
    \end{array}
   \]
  with  $0< \zeta_i <1$.    \\
  Using the same arguments as in the proof of Proposition \ref{etape4}, we obtain:
  \[
   \eE[Z_m(\efo+n^{-1/2} \uu,\efo) ] 
   =  \frac{f(0)}{2n} \sum^m_{i=1} [\uu^t \eg(\XX_i,\efo)]^2(1+o(1))=O(1).
  \]
   On the other hand, using  similar arguments to those for $W_m(\tau; \efo+n^{-1/2} \uu,\efo)$, we have that
   \[
   \Var [Z_m(\efo+n^{-1/2} \uu,\efo) ]
    \leq   \sum^m_{i=1} h^2_i(\efo+n^{-1/2} \uu,\efo) 
    =O(1),
       \]
   uniformly in $m$ and $\uu$.\\
  Then,  the Bienaym\'e-Tchebychev inequality, we obtain that $Z_m(\efo+n^{-1/2} \uu,\efo)=O_{\eP}(1)$. Taking into account  (\ref{ibis}) and  (\ref{Wmm}), the Lemma follows.
 \hspace*{\fill}$\blacksquare$ \\

\section{Quantile regression with multiple change-points}
\label{sec3:1}
\hh This section considers that the nonlinear model changes to unknown observations. More specifically, the regression parameters change to unknown times. First, we define the quantile estimators of the model parameters. If the number of changes is known, we give the convergence rate and the limiting distribution of the all estimators. Next, we give a consistent criterion for estimating the change-point number. \\

 Consider  a model with $K$ change-points, i.e. a model which  changes to observations $l_1, \cdots, l_K$, with $1 < l_1 < \cdots, l_K<n$,
\begin{equation}
\label{eth3}
Y_i= \sum^{K}_{r=0} g(\XX_i, \ef_{r+1}) \e1_{l_r \leq i < l_{r+1}}+ \varepsilon_i, 
\end{equation}
$i=1, \cdots , n$, with $l_0=1$ and $l_{K+1}= n$.\\
 We assume that  numbers of changes $K$ is known. \\

 Concerning the change-point location, we suppose that each segment contains a significant proportion of samples:\\
\textbf{(A8)} $l_{r+1} - l_r \geq n^a$, $a > 1/2$, for all $r=0, \cdots, K$, with $l_0=1$ and $l_{K+1}=n$.\\
This condition is necessary in order to apply Lemma \ref{Lemma2.3_adapt}, therefore constant "$a$" must be strictly greater than 1/2.\\
 
 For fixed $K$, the parameters of model (\ref{eth3}) are the regression parameters  $\eth_1 \equiv (\ef_1, \cdots, \ef_{K+1}) \in \Gamma^{K+1}$ and the change-points $\eth_2 \equiv (l_1, \cdots, l_K) \in \N^K$. The true values of the parameters are $\eth_1^0 \equiv (\ef^0_1, \cdots , \ef^0_{K+1})$ for the regression parameters and  $\eth_2^0 \equiv (l^0_1, \cdots , l^0_K)$ for the change-points. Obviously $\ef^0_{r+1} \neq \ef^0_r$, for all $r=1, \cdots , K$.\\
 
 We define the  quantile estimators of   parameters $\eth_1$ and $\eth_2$  by
\begin{equation}
\label{est_param_chp}
(\hat \eth^{(\tau)}_{1n}, \hat \eth^{(\tau)}_{2n}) \equiv 
 \argmin_{(\eth_1, \eth_2)}  \sum^{K+1}_{r=1} \sum^{l_r}_{i=l_{r-1}+1} \rho_\tau(Y_i - g(\XX_i,\ef_r)).
\end{equation}
See \cite{Ciuperca:11} for a discussion on the construction of the estimators in a change-point model. 
\\

In order to prove the convergence rate of the change-point quantile estimator $\hat \eth^{(\tau)}_{2n}$, we first prove that if in a phase we take in the place of the true regression parameter those of the nearby phase, then the value of the objective function is different from that calculated for the true value. 

\begin{lemma}
\label{Lemma3.1_adapt}
Under assumption (A6), we have for every $r=1, \cdots, K$, when $l_r < l^0_r$ such that $l_r^0-l_r \rightarrow\infty$, that there exists $\eta >0$, $C>0$ such that
\[
 \eP\left[\left|\sum^{l^0_r}_{i=l_r+1}(\rho_\tau(\varepsilon_i- g(\XX_i,\ef^0_{r+1})+g(\XX_i,\ef^0_r))-\rho(\varepsilon_i)) \right|
  \geq \eta(l^0_r-l_r)  \right] \geq 1 - \exp(-C (l^0_r-l_r)).
\]
\end{lemma}
{\bf Proof.}
Since $\ef^0_{r+1} \neq \ef^0_r$, there exists  $\delta$, and $ \epsilon_0 >0$ such that $|g(\XX_i,\ef^0_{r+1})-g(\XX_i,\ef^0_r)| \geq \delta$ for $(l^0_r-l_r)\epsilon_0$ observations. Then
$$
\eE\left[ \sum^{l^0_r}_{i=l_r+1} G^{(\tau)}_i(\ef^0_{r+1},\ef^0_r)\right]=(l^0_r-l_r)\epsilon_0 \int^0_{- \delta }(x+\delta) dF(x).$$
Applying Proposition \ref{Lemma2.4_adapt}(i) for $c=\| \ef^0_{r+1}- \ef^0_r \|_2$ and  $a_n=l^0_r-l_r$, we have that for all  $\epsilon >0 $, the following inequality
$$
 \eP \left[\left|\sum^{l^0_r}_{i=l_r+1} G^{(\tau)}_i(\ef^0_{r+1},\ef^0_r)- \eE[ \sum^{l^0_r}_{i=l_r+1} G^{(\tau)}_i(\ef^0_{r+1},\ef^0_r)]\right| 
 \geq \epsilon(l^0_r - l_r)\right] \leq \exp(- C \epsilon^2 (l^0_r-l_r)).
$$
 Then, Lemma follows  considering  $\eta=\epsilon=2^{-1} \epsilon_0 \int^0_{- \delta }(x+\delta) dF(x)$.
 \hspace*{\fill}$\blacksquare$ \\

\begin{remark}
\label{Remark 2.11_adapt}
Using Lemma \ref{Lemma 5_LAD} and Proposition \ref{Lemma2.4_adapt}, by  similar technique to one use in the paper of \cite{Ciuperca:11}, for Lemmas 7 and 8, and in the paper \cite{Ciuperca:14}, for Lemmas 3 and 4, we obtain their equivalent. That is, if data come from two different models, the quantile estimator is close to the parameter of the model from where most of the data came. 

\end{remark}

Following result shows that the distance between the change-point quantile  estimator and the true value is finished.
\begin{theorem}
\label{Theorem3.1_adapt} Under assumptions (A1), (A4), (A6)-(A8) and  if density function $f$ of $\varepsilon$ satisfies conditions of Proposition \ref{etape4}, then we have 
$\| \hat \eth^{(\tau)}_{2n}-\eth^0_2\|_2=O_{\eP}(1)$.
\end{theorem}
{\bf Proof.}
The proof is similarly to that of Theorem 3.1 of  \cite{Ciuperca:14}, using   relation (\ref{eq38_adapt}), $\eE[G_i^{(\tau)}(\ef,\efo) ] \geq 0$ by (\ref{egn}) and   Lemma \ref{Lemma2.3_adapt}, Proposition \ref{etape4}, Lemma \ref{Lemma3.1_adapt}, Remark \ref{Remark 2.11_adapt}. We omit all details. 
\hspace*{\fill}$\blacksquare$ \\

With this result we can now give the asymptotic distributions, first for the change-point estimator and then for the regression parameter estimator. This result is the generalization of that obtained in \cite{Ciuperca:11} for LAD method ($\tau=1/2$), where the proof was based on norm $L_1$ of objective function. The asymptotic distribution of the change-point quantile estimator depends on regression function $g$, on the true regression parameters to the left and right of the estimated break point and of quantile index $\tau$ of $\varepsilon$. The asymptotic distribution of regression parameter quantile estimator  is Gaussian, with covariance matrix dependent of $\tau$. Theorem \ref{Theorem2_Oka_Qu} is a standard result in quantile model without change-point (\cite{Koenker:05}) and in a change-point model estimated by other methods. See i.e. \cite{Boldea:Hall:13} for LS method.   \\

We consider by convention that $l^0_0= \hat l^{(\tau)}_0=1$ and $l^0_{K+1}=\hat l^{(\tau)}_{K+1} =n$. 

\begin{theorem}
\label{Theorem2_Oka_Qu}
Under the same  conditions of Theorem \ref{Theorem3.1_adapt},  we have the following asymptotic laws of the change-point quantile estimators: \\
 (i) for each $r=1, \cdots, K$, $$ (\hat l^{(\tau)}_r- l^0_r)  \overset{{\cal L}} {\underset{n \rightarrow \infty}{\longrightarrow}} \argmin_{j \in \Z} Z^{(\tau)}_{r,j},$$
  where:\\
- if $j=1, 2, \cdots $,
\[
Z^{(\tau)}_{r,j} \equiv \sum^{l^0_r+j}_{i=l^0_r+1} \left[\rho_\tau(\varepsilon_i -g(\XX_i,\ef^0_r)   +g(\XX_i,\ef^0_{r+1}))-\rho_\tau(\varepsilon_i)\right].
\]
- if $j= -1, -2, \cdots, $ 
\[
Z^{(\tau)}_{r,j}\equiv \sum^{l^0_r}_{i=l^0_r+j} \left[\rho_\tau(\varepsilon_i -g(\XX_i,\ef^0_{r+1})
 +g(\XX_i,\ef^0_{r}))-\rho_\tau(\varepsilon_i)\right].
\]
(ii) for each $r=1, \cdots, K+1$, 
$$(\hat l^{(\tau)}_r -\hat l^{(\tau)}_{r-1})^{1/2}(\hat \ef^{(\tau)}_r - \ef^0_{r}) \eS^{1/2}_r  \overset{{\cal L}} {\underset{n \rightarrow \infty}{\longrightarrow}} {\cal N} (\textbf{0}, \frac{\tau (1-\tau)}{f^2(0)}\textbf{I}_p ),$$
 with $$\eS_r \equiv (l^0_r-l^0_{r-1})^{-1} \sum^{l^0_r}_{i={l^0_{r-1}+1}} \eg(\XX_i,\ef^0_r)\eg^t(\XX_i,\ef^0_r)$$ and $\textbf{I}_p$ the identity matrix of order $p$.
\end{theorem}
 {\bf Proof.} 
Let us consider the set of change-point vectors 
 \[
 \Theta_2 \equiv \{ \eth_2=(l_1, \cdots, l_K); 
   l_j=l^0_j+m, |m| \leq C_2 , \forall j=1, \cdots , K  \}
 \]
  and the set of regression parameter vectors
   \[
   \Theta_1 \equiv \{\eth_1=(\ef_1, \cdots \ef_K) ; 
  (l^0_j-l^0_{j-1})^{1/2}\| \ef_j-\ef^0_j \|_2 \leq C_3, \forall j=1, \cdots , K  +1 \},
      \]
   with $C_2, C_3 >0$ finite constants. Let be the sum 
   $$S_n(\tau,\eth_1, \eth_2) \equiv \sum^K_{j=0} \sum^{l_{j+1}}_{i=l_j+1}  \rho_\tau(Y_i-g(\XX_i, \ef_{j+1})) .$$
    Consider the following identity
\begin{equation}
 \label{et1}
  \inf_{\eth_2} \inf_{\eth_1} S_n(\tau,\eth_1, \eth_2) = \inf_{\eth_2}\inf_{\eth_1} \pth{ S_n(\tau,\eth_1, \eth_2^0)
  +  S_n(\tau,\eth_1, \eth_2)-  S_n(\tau,\eth_1, \eth_2^0) } .
 \end{equation}
By Lemma \ref{Lemma L+} we have that $[S_n(\tau,\eth_1, \eth_2)-  S_n(\tau,\eth_1, \eth_2^0) ]- [S_n(\tau,\eth^0_1, \eth_2)-  S_n(\tau,\eth^0_1, \eth_2^0)] = o_{\eP}(1)$ uniformly in  $\eth_1,\eth_2$ belonging in $ \Theta_1 \times \Theta_2$. \\
\hh Without loss of generality, we suppose that $ \hat l^{(\tau)}_r \leq l^0_r$.\\
 By the definition of $S_n$,  we have that 
 $$
  S_n(\tau,\eth^0_1, \eth_2)-  S_n(\tau,\eth^0_1, \eth_2^0) =  \sum^K_{r=1} \sum^{l^0_r}_{i=l_r+1}
  \left[\rho_\tau(Y_i-g(\XX_i, \ef^0_{r+1}))-\rho_\tau(Y_i -g(\XX_i, \ef^0_r))\right] .
  $$
Then, relation (\ref{et1}) becomes
$$
 \inf_{\eth_1 \in \Theta_1} S_n(\tau,\eth_1,\eth_2^0) +  \inf_{\eth_2 \in \Theta_2} \sum^K_{r=1}  \sum^{l^0_r}_{i=l_r+1} \left[\rho_\tau(Y_i-
 g(\XX_i, \ef^0_{r+1})) -\rho_\tau(Y_i -g(\XX_i, \ef^0_r))\right](1+o_{\eP}(1)).
$$
Theorem results taking into account that  every term of this last  relation  depends on different parameters, together with convergence rate of the  estimators (by Theorem \ref{Theorem3.1_adapt} for the change-point estimator and $(l^0_r-l^0_{r-1})^{-1/2}$ for the regression parameter quantile estimator) and limit law of quantile estimator for a nonlinear model  (see for example \cite{Koenker:05}).
\hspace*{\fill}$\blacksquare$ \\

\begin{remark}
 In the case presented here, parameters  $\ef_r$, $\ef_{r+1}$ from a segment to the other are fixed. In the paper of \cite{Oka:Qu:11} for linear model, it is supposed that the difference between two consecutive parameter tends to zero as $n \rightarrow \infty$. Then, the limit law of the change-points estimators  is   totally different, it is the maximizer of a Wiener process with drift.
\end{remark}
\begin{remark}
\label{remark5}
In order to determine the number of change points, we can use a similar criterion to that proposed in the paper of \cite{Ciuperca:14} for a linear quantile model.   Under conditions that  $\eE[\rho_\tau(\varepsilon)]>0$ and $\eE[\rho^2_\tau(\varepsilon)] < \infty$, we propose the following  consistent  estimator of the change-point number $K$ 
 \begin{equation}
 \label{eth4}
  \hat K^{(\tau)}_n \equiv  
\argmin_K \left( n \log \left(n^{-1} S_n(\tau,\hat \eth^{(\tau)}_{1n}(K),\hat \eth^{(\tau)}_{2n}(K)) \right)
 +P(K,p)B_n \right), 
 \end{equation}
where the function $S_n$ is defined in the proof of Theorem \ref{Theorem2_Oka_Qu}, $(\hat \eth^{(\tau)}_{1n}(K),\hat \eth^{(\tau)}_{2n}(K))$ is  the quantile estimators of $(\eth_1,\eth_2)$ for a fixed $K$, $(B_n)$ is  a deterministic sequence converging to infinity such that $B_n n^{-a} \rightarrow 0$, $B_n n^{-1/2} \rightarrow \infty$, as $n \rightarrow \infty$ and the penalty function $P(K,p)$ is such that $P(K_1,p) \leq P(K_2,p)$ for all number change-points $K_1 \leq K_2$. Recall that the constant $a$ is that of the supposition (A8) and $p$ is parameter number of the regression function $g$. \\
The proof of the consistency of the criterion  is similar to that in \cite{Ciuperca:14}. We do not give the details.
\end{remark}
\section{Simulation study for change-point nonlinear models}
\label{sec4:1} 
\hh To evaluate the performance of the quantile method in a change-point nonlinear  model,  Monte Carlo simulations  are realized. We compare the performance of the least squares (LS)  and quantile estimation methods. We use\textit{ quantreg, VGAM packages in R} to run the simulations.\\

 For each model, $100$ Monte Carlo samples of size $n$ are generated for regressor $X$ and error $\varepsilon$. \\
 Throughout this section,  we generate the design $X \sim {\cal N}(1,1)$ and the regression function $g(x, \ef)$ is growth function $b_1-\exp(-b_2x)$, or more exactly the mono-molecular model (see \cite{Seber:Wild:03}), with $\ef=(b_1,b_2)$. The same regression function has been considered in \cite{Ciuperca:11a} using the M-method that has the least squares method as a special case.  For the errors $\varepsilon$, three distributions were considered: standard Normal ${\cal N}(0,1)$, Laplace ${\cal L}(0,1)$, and Cauchy ${\cal C}(0,1)$.   \\
  
The quantile estimations of the regression parameters and of the change-points, for a fixed number $K$ of change-points, are calculated using relation  (\ref{est_param_chp}). The corresponding LS estimations are obtained by minimizing in $\eth_1$ and $\eth_2$ following sum (see \cite{Boldea:Hall:13}):
$$\sum^K_{r=1} \sum^{l_r}_{i=l_{r-1}+1} \left[Y_i- g(\XX_i,\ef_r)\right]^2.$$
 
 \subsection{Known change-point number}
\hh First, in Tables \ref{Tabl1}, \ref{Tabl2}, \ref{Tabl3}, the change-point number is known and it is equal to two (model with three phases). In Tables \ref{Tabl1} and \ref{Tabl2} the number of observations is n=100, with the particular case of epidemic model in Tables \ref{Tabl2},   when model is the same in the first and the third phase ($\ef^0_1=\ef^0_3=(0.5, 1)$).\\
 Since the asymptotic distribution of the change-point quantile estimators can not be symmetric (see Theorem \ref{Theorem2_Oka_Qu}), the median of change-point estimations are given. Asymptotic distributions of regression parameter estimators by LS and quantile methods in a change-point nonlinear model are Gaussian (see Theorem \ref{Theorem2_Oka_Qu} and corresponding result of \cite{Boldea:Hall:13} for LS method). Then, the  mean and standard-deviation(\textit{sd)} of corresponding estimations are reported.\\
  In all situations (see Tables \ref{Tabl1}, \ref{Tabl2}, \ref{Tabl3}), the median of the change-point estimations are very close to the true values. 
When the errors are Gaussian , very  good results are obtained by the two estimation methods. For  Laplace errors, the results deteriorate slightly, while for Cauchy errors, the quantile method gives very satisfactory results, while by LS method, the obtained estimates are biased and with a wide variation, when $n=100$ or when $n$ is greater (Tables \ref{Tabl1} and \ref{Tabl3}). 
 
\begin{table*}[t]
{\scriptsize
\caption{\footnotesize Model with two change-points $l^0_1=20$, $l^0_2=85$, $n=100$.  $\ef^0_1=(0.5, 1)$, $ \ef^0_2=(1, -0.5)$, $ \ef^0_3=(2.5, 1)$. Estimation by LS and quantile methods.}
\begin{center}
\begin{tabular}{|c|c|c|c|c|c|c|}\hline
Estimation & $\varepsilon$ law & median($\hat l^{(\tau)}_1$) & median($\hat l^{(\tau)}_2$) & mean($\hat \ef^{(\tau)}_1$) & mean($\hat \ef^{(\tau)}_2$) & mean($\hat \ef^{(\tau)}_3$) \\
method & & & &  \textit{sd}($\hat \ef^{(\tau)}_1$) & \textit{sd}($\hat \ef^{(\tau)}_2$) & \textit{sd}($\hat \ef^{(\tau)}_3$) \\ \hline
 LS & $\varepsilon \sim {\cal N}$ & 19 & 84 & (0.52, 1.06) & (0.98, -0.5) & (2.52, 1.07) \\
  & & & & \textit{(0.17, 0.64)} & \textit{(0.09, 0.02)} & \textit{(0.15, 0.5)} \\
\cline{2-7} 
   & $\varepsilon \sim {\cal L}$ & 19 & 84 & (0.58, 1.28) & (0.98, -0.5) & (2.65, 1.13) \\
  & & & & \textit{(0.42, 1.56) }& \textit{(0.22, 0.05)} & \textit{(0.46, 1.1)} \\
  \cline{2-7}
  & $\varepsilon \sim {\cal C}$ & 22 & 85 & (2.51, 1.26) & (2.34, -0.24) & (7.7, 1.75) \\
  & & & & \textit{(18.7, 2.2)} & \textit{(12.7, 0.96)} & \textit{(42, 3.4)} \\ \hline
   quantile & $\varepsilon \sim {\cal N}$ & 19 & 84 & (0.52, 1.1) & (0.99, -0.5) & (2.53, 1.1) \\
  & & & & \textit{(0.16, 0.78) }& \textit{(0.09, 0.02)} & \textit{(0.18, 0.8)} \\
\cline{2-7} 
   & $\varepsilon \sim {\cal L}$ & 19 & 84 & (0.57, 1.17) & (1, -0.5) & (2.6, 1.45) \\
  & & & & \textit{(0.37, 1.28)} & \textit{(0.13, 0.04)} & \textit{(0.32, 3.2)} \\
  \cline{2-7}
  & $\varepsilon \sim {\cal C}$ & 20 & 84 & (0.58, 1.2) & (0.98, -0.48) & (2.7, 1.75) \\
  & & & & \textit{(0.55, 1.1) }& \textit{(0.29, 0.23)} & \textit{(0.6, 3.1)} \\ \hline
\end{tabular}
 \end{center}
\label{Tabl1} 
}
\end{table*}
\begin{table*}[t] 
{\scriptsize
\caption{\footnotesize Model with two change-points $l^0_1=20$, $l^0_2=85$, $n=100$.  $\ef^0_1=\ef^0_3=(0.5, 1)$, $ \ef^0_2=(1, -0.5)$. Estimation by LS and quantile methods.}
\begin{center}
\begin{tabular}{|c|c|c|c|c|c|c|}\hline
Estimation & $\varepsilon$ law & median($\hat l^{(\tau)}_1$) & median($\hat l^{(\tau)}_2$) & mean($\hat \ef^{(\tau)}_1$) & mean($\hat \ef^{(\tau)}_2$) & mean($\hat \ef^{(\tau)}_3$) \\
method & & & &  \textit{sd}($\hat \ef^{(\tau)}_1$) & \textit{sd}($\hat \ef^{(\tau)}_2$) & \textit{sd}($\hat \ef^{(\tau)}_3$) \\ \hline
 LS & $\varepsilon \sim {\cal N}$ & 19 & 84 & (0.5, 1.02) & (1, -0.5) & (0.5, 1.2) \\
  & & & & \textit{(0.12, 0.31)} & \textit{(0.07, 0.02)} & \textit{(0.14, 1.7)} \\
\cline{2-7} 
   & $\varepsilon \sim {\cal L}$ & 19 & 84 & (0.58, 1.03) & (1.03, -0.51) & (0.5, 2.01) \\
  & & & & \textit{(0.35, 0.62)} & \textit{(0.24, 0.05)} & \textit{(0.39, 6.5) }\\
  \cline{2-7}
  & $\varepsilon \sim {\cal C}$ & 21 & 85 & (2.28, 1.75) & (-0.07, -0.07) & (1.3, 1.8) \\
  & & & & \textit{(11.8, 3)} & \textit{(17, 1.2)} & \textit{(7.9, 6.4)} \\ \hline
   quantile & $\varepsilon \sim {\cal N}$ & 19 & 84 & (0.52, 1.05) & (1, -0.5) & (0.5, 1.02) \\
  & & & & \textit{(0.15, 0.52) }& \textit{(0.09, 0.02)} & \textit{(0.17, 0.4)} \\
\cline{2-7} 
   & $\varepsilon \sim {\cal L}$ & 19 & 84 & (0.56, 1.06) & (1.03, -0.5) & (0.51, 1.8) \\
  & & & & \textit{(0.3, 0.8)} & \textit{(0.15, 0.04)} & \textit{(0.3, 5.6)} \\
  \cline{2-7}
  & $\varepsilon \sim {\cal C}$ & 19 & 84 & (0.64, 1.45) & (0.93, -0.49) & (0.7, 1.3) \\
  & & & & \textit{(0.56, 1.6) }& \textit{(0.2, 0.14)} & \textit{(0.56, 1.54)} \\ \hline
\end{tabular}
 \end{center}
\label{Tabl2} 
}
\end{table*}

\begin{table*}[t] 
{\scriptsize
\caption{\footnotesize Model with two change-points $l^0_1=100$, $l^0_2=200$, $n=300$.  $\ef^0_1=(0.5, 1)$, $ \ef^0_2=(1, -0.5)$, $ \ef^0_3=(2.5, 1)$. Estimation by LS and quantile methods.}
\begin{center}
\begin{tabular}{|c|c|c|c|c|c|c|}\hline
Estimation & $\varepsilon$ law & median($\hat l^{(\tau)}_1$) & median($\hat l^{(\tau)}_2$) & mean($\hat \ef^{(\tau)}_1$) & mean($\hat \ef^{(\tau)}_2$) & mean($\hat \ef^{(\tau)}_3$) \\
method & & & &  \textit{sd}($\hat \ef^{(\tau)}_1$) & \textit{sd}($\hat \ef^{(\tau)}_2$) & \textit{sd}($\hat \ef^{(\tau)}_3$) \\ \hline
 LS & $\varepsilon \sim {\cal N}$ & 99 & 199 & (0.5, 1) & (1, -0.5) & (2.5, 1) \\
  & & & & \textit{(0.05, 0.07)} & \textit{(0.06, 0.01)} & \textit{(0.05, 0.08)} \\
\cline{2-7} 
   & $\varepsilon \sim {\cal L}$ & 99 & 199 & (0.48, 1.02) & (1, -0.5) & (2.5, 1) \\
  & & & & \textit{(0.16, 0.19) }& \textit{(0.17, 0.05)} & \textit{(0.13, 0.16)} \\
  \cline{2-7}
  & $\varepsilon \sim {\cal C}$ & 100 & 200 & (3.05, 1.02) & (1.99, -0.32) & (3.2, 1.1) \\
  & & & & \textit{(16.2, 1.17} & \textit{(10.6, 0.77)} & \textit{(4.9, 1.06)} \\ \hline
   quantile & $\varepsilon \sim {\cal N}$ & 99 & 199 & (0.5, 1) & (1, -0.5) & (2.5, 1) \\
  & & & & \textit{(0.05, 0.08)} & \textit{(0.07, 0.02)} & \textit{(0.05, 0.08)} \\
\cline{2-7} 
   & $\varepsilon \sim {\cal L}$ & 99 & 199 & (0.51, 1) & (0.99, -0.5) & (2.5, 1.01) \\
  & & & & \textit{(0.09, 0.16) }& \textit{(0.11, 0.03)} & \textit{(0.1, 0.1)} \\
  \cline{2-7}
  & $\varepsilon \sim {\cal C}$ & 99 & 199 & (0.54, 1.06) & (0.98, -0.5) & (2.5, 0.98) \\
  & & & & \textit{(0.14, 0.75)} & \textit{(0.16, 0.04)} & \textit{(0.1, 0.4)} \\ \hline
\end{tabular}
 \end{center}
\label{Tabl3} 
}
\end{table*}

\begin{table*}[t] 
{\scriptsize
\caption{\footnotesize Results on the choice of the change-point number by criteria associated to  LS and quantile methods. The true change-point number is 1 in $l^0_1=20$ for $n=100$ observations. 100 Monte Carlo replications.}
\begin{center}
\begin{tabular}{|c|c|c|c|c|c|c|c|c|c|c|c|c|}\hline
True parameters &  \multicolumn{6}{c|} {LS method} &  \multicolumn{6}{c|} {Quantile method} \\
\cline{2-13}
 & \multicolumn{3}{c|} {$\varepsilon \sim {\cal N}$} & \multicolumn{3}{c|} {$\varepsilon \sim {\cal C}$} & \multicolumn{3}{c|} {$\varepsilon \sim {\cal N}$} & \multicolumn{3}{c|} {$\varepsilon \sim {\cal C}$} \\ 
  & $\hat K^{(\tau)}_n=0$ & $=1$ & $=2$  & $=0$ & $=1$ & $=2$  & $=0$ & $=1$ & $=2$ & $=0$ & $=1$ & $=2$ \\ \hline
  $\ef^0_1=(0.5, 1)$, $\ef^0_2=(10, 2.5)$ & 0& 100 & 0 & 56 &43 &1 & 2 & 92 & 6 & 13 & 86 & 1 \\ \hline
  $\ef^0_1=(0.5, 1)$, $\ef^0_2=(1, -0.5)$ & 4& 96 & 0 & 91 &9 & 0 &5 & 95 & 0 & 92 & 8 & 0 \\ \hline
\end{tabular}
 \end{center}
\label{Tabl4} 
}
\end{table*}

\begin{figure*}[t]
\begin{minipage}[b]{0.45\linewidth}
 \centering \includegraphics[scale=0.45]{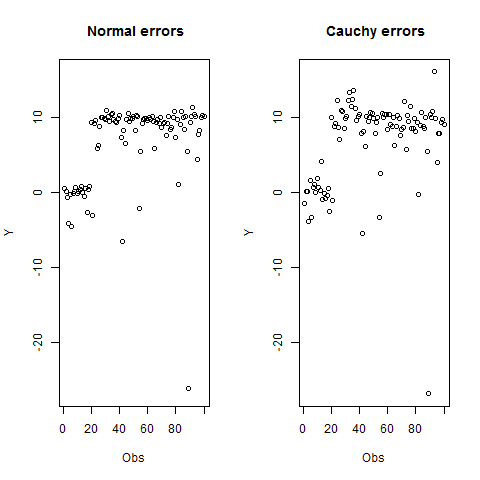}

\caption{\it   Model with normal and Cauchy errors, with one change-point in $l^0=20$, far  parameters for the two phases.}
\label{Figure 1}
\end{minipage} \hfill
   \begin{minipage}[b]{0.45\linewidth}
\centering \includegraphics[scale=0.45]{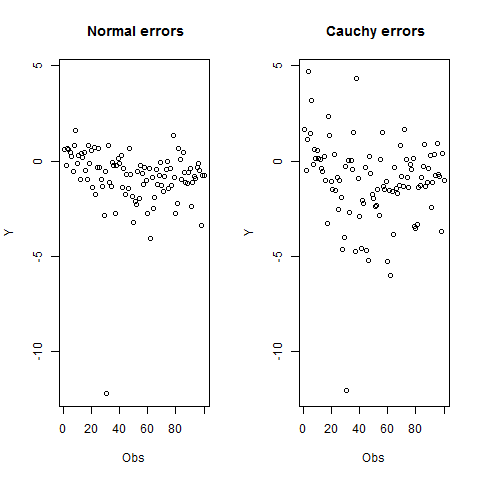}
\caption{ \it  Model with normal and Cauchy errors, with one change-point in $l^0=20$, close parameters for the two phases.}
\label{Figure 2}  
   \end{minipage}
\end{figure*}

\subsection{Unknown change-point number}
\hh In view of these results, in order to study the selection criterion of the change-point  number, we will consider only Normal and Cauchy distributions for errors. We simulate a model with one change-point in $l^0_1=20$ for $n=100$ observations. The  estimation $\hat K^{(\tau)}_n$ of the change-point number associated to quantile method  is calculated using Remark \ref{remark5}. For criterion (\ref{eth4}), for the penalty we consider  $P(K,p)=Kp$ and deterministic sequence $B_n=n^{5/8}$.  The estimation of the change-point number associated to LS method, is the minimizer in $K$ of
\[
 n \log \left(n^{-1} \min_{(\eth_1,\eth_2)}\acc{ \sum^K_{r=1} \sum^{l_r}_{i=l_{r-1}+1} \left[Y_i- g(\XX_i,\ef_r)\right]^2} \right)
 +P(K,p)B_n.
\]
\hh   Two cases are considered for the true regression parameters: the parameters of the two (true) phases are far (Figure \ref{Figure 1}), $\ef^0_1=(0.5, 1)$, $\ef^0_2=(10, 2.5)$ and the parameters are closely (Figure \ref{Figure 2}),  $\ef^0_1=(0.5, 1)$, $\ef^0_2=(1, -0.5)$. In the case of Gaussian errors, the criterion associated to the  LS method  is slightly better when the parameters are far. The criteria associated to the two methods (LS and quantile) give the same good results if the parameters are closely. In the case of Cauchy errors, the quantile criterion selects well the change-point number when the parameters are far, while when the parameters are closely, the two criteria rather prefer a model without change-points (Tables \ref{Tabl4}).
\subsection{Conclusion}
\hh These simulations  allow us to conclude that for a nonlinear model with change-points, when the errors are Gaussian, the quantile method, proposed in this paper, gives similar results to those obtained by least squares method. On the other hand, for heavy-tailed errors, the performance of the  quantile method is better than LS method, whether in estimation or in selection criterion.  
  \appendix
\section{Bernstein's Inequality}
\textbf{Bernstein's Inequality} (see for example \cite{Pollard:84}).\\
Let $Z_i$ be a sequence of independent random variables with mean zero and $|Z_i| \leq \beta $ for some $\beta >0$. Let also $V \geq \sum^n_{i=1}\eE[Z_i^2] $. Then for all $0<s<1$ and $0 \leq z \leq V/(s \beta)$, we have
\begin{equation}
\label{ineq_Bernstein}
\eP\cro{| \sum^n_{i=1}Z_i| >z} \leq 2 \exp \pth{-z^2s(1-s)/V}.
\end{equation}


\begin{thebibliography}{3}
\bibitem[Bai(1998)]{Bai:98}
Bai, J., 
\newblock Estimation of multiple-regime regressions with least absolute deviation.
\newblock {\em Journal of Statistical Planning Inference}, {\bf 74},  103-134, (1998).
\bibitem[Bai and Perron(1998)]{Bai:Perron:98} 
Bai,  J.,  Perron P., 
\newblock Estimating and testing linear models with multiple structural changes,
\newblock {\em Econometrica} {\bf 66}(1), 47-78, (1998).
\bibitem[Boldea and Hall(2013)] {Boldea:Hall:13}
Boldea, O., Hall, A.R.,  
\newblock {Estimation and inference in unstable nonlinear least squares models}. 
 \newblock {\textit{Journal of Econometrics}},   \textbf{172}(1), 158-167, (2013). 
\bibitem[Chen et al.(2013)]{Chen:Tran:Lin:04}
Chen L.A., Tran L.T., Lin L.C., 
\newblock Symmetric regression quantile and its application to robust estimation for the nonlinear regression model, 
\newblock \textit{Journal of Statistical Planning and Inference},  \textbf{126}(2), 423-440, (2004).
\bibitem[Choi et al.(2005)]{Choi:Kim:Lee:05}
Choi S.H., Kim K.J., Lee M.S., 
\newblock Robust test based on nonlinear regression quantile estimators, 
\newblock \textit{Communications of the Korean Mathematical Society},  \textbf{20}(1), 145-159, (2005).
\bibitem[Ciuperca(2009)]{Ciuperca:09}
Ciuperca G., 
\newblock The  M-estimation in a multi-phase random nonlinear model. 
\newblock {\em Statistics and  Probability Letters}, {\bf 75}(5), 573-580, (2009).
\bibitem[Ciuperca(2011a)]{Ciuperca:11a}
Ciuperca, G., 
\newblock  A general criterion to determinate the number of change-points. 
\newblock  \textit{Statistics and Probability Letters},
 \textbf{81}, no 8, 1267-1275, (2011).
\bibitem[Ciuperca(2011b)]{Ciuperca:11}
Ciuperca G., 
\newblock Estimating nonlinear regression with and without change-points by the LAD-method.
\newblock{\it Annals of the Institute of Statistical Mathematics},   \textbf{63}(4), 717-743, (2011).
\bibitem[Ciuperca(2014)]{Ciuperca:14}
Ciuperca G.,  
\newblock Adaptive model selection in a high-dimension multiphase quantile regression.
\newblock{\it arXiv:1309:1262}, (2014).
\bibitem[Furno(2012)]{Furno:12}
Furno M., 
\newblock Tests for structural break in quantile regressions.
\newblock{\it AStA Advances in Statistical Analysis},   \textbf{96}(4), 493-515, (2012).
\bibitem[Knight(1998)]{Knight:98}
Knight K., 
\newblock {Limiting distributions for $L_1$ regression estimators under general conditions}, 
\newblock {\it Annals of Statistics}, \textbf{26}(2), 755-770, (1998).
\bibitem[Koenker(2005)]{Koenker:05}
Koenker R.,
\newblock {Quantiles regression}, 
\newblock {Econometric Society Monographs}, No \textbf{38}, Cambridge University Press, (2005).
\bibitem[Koenker and Bassett(1978)]{Koenker:Bassett:78}
Koenker R., Bassett G., 
\newblock {Regression Quantiles}, 
\newblock {\it Econometrica}, \textbf{46}, 33-50, (1978).
\bibitem[Koul and Qian(2002)]{Koul:Qian:02}
Koul, H.L.,   Qian, L., 
\newblock  Asymptotics of maximum likelihood estimator in a two-phase linear regression model.
\newblock {\em Journal of Statitical Planning and Inference}  {\bf 108}, 99-119, (2002).
\bibitem[Koul et al.(2003)]{KQS}
Koul, H.L.,  Qian, L.,  Surgailis, D., 
\newblock  Asymptotics of M-estimators in two-phase linear regression models.
\newblock {\em Stochastic Processes and their Applications } {\bf 103}, 123-154, (2003).
\bibitem[Oberhofer(1982)]{Oberhofer:82}
Oberhofer W.,  
\newblock The consistency of nonlinear regression minimizing the L1-norm.
\newblock {\it Annals of Statistics}, \textbf{10}(1),  316-319, (1982).
\bibitem[Oberhofer and Haupt(2014)]{Oberhofer:Haupt:13}
 Oberhofer W., Haupt H.,  
\newblock Asymptotic theory for nonlinear quantile regression under weak dependence.
\newblock Working Paper. \textit{https://www.researchgate.net/publication/29858831\_Asymptotic\_theory\_for\_nonlinear}\\
\textit{\_quantile\_regression\_under\_weak\_dependence}, (2014).
\bibitem[Oka and Qu(2011)]{Oka:Qu:11}
Oka  T.,  Qu Z.,  
\newblock Estimating structural changes in regression quantiles.
\newblock {\it Journal of Econometrics}, \textbf{162},  248-267, (2011).
\bibitem[Pollard(1984)]{Pollard:84}
Pollard D.,  
\newblock \textit{Convergence of stochastic processes}, 
\newblock{Springer}, New York, (1984).
\bibitem[Qu(2008)]{Qu:08}
 Qu Z.,  
\newblock Testing for structural change in regression quantiles.
\newblock {\it Journal of Econometrics}, \textbf{146},  170-184, (2008).
\bibitem[Seber and Wild(2003)]{Seber:Wild:03}
 Seber,  G.A.F.,  Wild, C.J., 
\newblock  {\em Nonlinear regression}.
\newblock {Wiley Series in Probability and Statistics}, New Jersey, (2003).
\bibitem[van der Vaart and Wellner(1996)]{van der Vaart:Wellner:96}
van der Vaart A.W., Wellner J.A.,. 
\newblock \textit{Weak convergence and empirical processes. With applications to statistics}, 
\newblock{Springer Series in Statistics}, New York, (1996).
\end{thebibliography}
\end{document}